\documentclass[12pt,a4paper, twoside, reqno]{amsart} 

\usepackage[margin=3cm]{geometry}

\usepackage{ifpdf}
\usepackage{graphicx}
\usepackage[T1]{fontenc} 
\usepackage[latin1]{inputenc} 
\usepackage{amsmath} 
\usepackage{amsthm} 
\usepackage{amsfonts} 
\usepackage{amssymb}
\usepackage{amsbsy}
\usepackage{bbm} 
\usepackage[english]{babel} 
\usepackage{varioref} 
\usepackage{enumerate} 

\theoremstyle{definition}

\theoremstyle{remark}
\newtheorem*{rem}{Remark}
\theoremstyle{plain}
\newtheorem{thm}{Theorem}
\newtheorem*{thm*}{Theorem}

\newtheorem{cor}{Corollary}

\newtheorem*{thmktv}{Theorem KTV}
\newtheorem*{lembl}{Lemma BL}

\newcommand{\norm}[1]{\ensuremath{\left\Vert #1 \right\Vert}}

\newcommand{\infabs}[1]{\ensuremath{\left\vert #1 \right\vert}}
\newcommand{\infeuc}[1]{\ensuremath{\left\vert #1 \right\vert_2}}

\newcommand{\hdim}{{\dim}}
\newcommand{\Bad}{\text{\textup{\bf Bad}}}

 \DeclareMathOperator{\mat}{Mat}

\title{On shrinking targets for $\mathbb{Z}^m$ actions on
  tori}

\author{Yann Bugeaud, \ Stephen Harrap, \ Simon Kristensen  \\ and  \  Sanju Velani}

\address{Y. Bugeaud, Université Louis Pasteur, Mathématiques, 7, rue
  René Descartes, F-67084 Strasbourg cedex, France}

\address{S. Harrap, Department of Mathematics, University of York,
  Heslington, York YO10 5DD, United Kingdom}

\address{S. Kristensen, Department of Mathematical Sciences, Faculty
  of Science, University of Aarhus, Ny Munkegade, Building 530,
  DK-8000 Aarhus C, Denmark}

\address{S. Velani, Department of Mathematics, University of York,
  Heslington, York YO10 5DD, United Kingdom}

\date{}

\thanks{SK is a Steno Research Fellow funded by the Danish Natural
  Science Research Council.}

\thanks{SV: Research supported by EPSRC grants EP/E061613/1 and
EP/F027028/1. \\   }

\thanks{Mathematics Subject Classification 2000:  11K60, 11J83, 37E10, 37E45.}

\begin{document}

\begin{abstract}
Let $A$ be an $n \times m$ matrix with real entries. Consider the
set $\Bad_A$ of $\mathbf{x} \in [0,1)^n$ for which there exists a
constant $c(\mathbf{x})>0$ such that for any $\mathbf{q} \in
\mathbb{Z}^m$ the distance between $\mathbf{x}$ and the point $\{A
\mathbf{q}\}$ is at least $c(\mathbf{x})
\infabs{\mathbf{q}}^{-m/n}\!\!$.  It is shown that the
intersection of $\Bad_A$ with any suitably regular fractal set is
of maximal Hausdorff dimension. The linear form systems
investigated in this paper are natural extensions of irrational
rotations of the circle. Even in the latter one-dimensional case,
the results obtained are new.
\end{abstract}

\maketitle

\section{Introduction}
\label{sec:introduction}

Consider initially a rotation of the  unit circle through an angle
$\alpha$. Identifying the circle with the unit interval $[0,1)$
and the base point of the iteration with the origin, we are
considering the numbers $0, \{\alpha\}, \{2\alpha\}, \dots$ where
$\{\, . \,\}$ denotes the fractional part. If $\alpha $ is
rational, the rotation is periodic. On the other hand, it is a
classical result of Weyl \cite{MR1511862} that any irrational
rotation of the circle is ergodic. In other words,
$\{q\alpha\}_{q\in \mathbb{N}}$ is equidistributed for irrational
$\alpha$.

Almost every orbit of an ergodic transformation visits any fixed
set of positive measure infinitely often. The `shrinking target
problem' introduced in \cite{HV}
 formulates the natural question
of what happens if the target set -- the set of positive measure
--  is allowed to shrink with time. For example and more
precisely, is there an optimal `shrinking rate' for which almost
every orbit visits the shrinking target infinitely often? In the
specific case of irrational rotations of the circle, the shrinking
target sets correspond to subintervals of $[0,1)$ whose lengths
decay according to some specified function $\psi$.  In other
words, the problem translates to considering inequalities of the
type
\begin{equation}
  \label{eq:1}
  \norm{q\alpha - x} < \psi(q),
\end{equation}
where $x \in [0,1)$ and $\norm{\, . \, }$ denotes the distance to
the nearest integer. The following statement dates back to
Khintchine \cite{0060.12203} and gives the `optimal' choice of
$\psi$ in the non-trivial case that $\alpha$ is irrational and  $
x \neq s \alpha + t $ for any integers $s$ and $t$. The inequality
\begin{equation}
  \label{eq:12}
  \norm{q\alpha - x} < \frac{C(\alpha)}{q}
\end{equation}
is satisfied for infinitely many  integers  $q $  with $C(\alpha)
:= \tfrac{1}{4}\sqrt{1-4 \lambda(\alpha)^2}$   -- the quantity
$\lambda(\alpha) : = \liminf_{q \rightarrow \infty} q
\norm{q\alpha}$ is the Markoff constant of $\alpha$. Note that
$\lambda(\alpha)$ is strictly positive whenever $\alpha$ is badly
approximable by rationals.  Thus, the above statement strengthens
a result of   Minkowski \cite{31.0213.02}; namely that
\eqref{eq:12} has infinitely many solutions with $C(\alpha) =
\tfrac{1}{4}$. In the trivial case that $\alpha$ is irrational and
$ x = s \alpha + t $ for some integers $s$ and $t$,
the classical theorem of Hurwitz 
implies that the inequality
\begin{equation}
  \label{eq:11}
  \norm{q\alpha - x} < \frac{1 + \epsilon}{\sqrt{5} q}  \qquad (\epsilon >0)
\end{equation}
is satisfied for infinitely many  integers  $q $.  Since
\eqref{eq:11} is weaker than \eqref{eq:12},  it follows that for
any irrational $\alpha$ and any $x$ the inequality  \eqref{eq:11}
has infinitely many solutions. We now describe a metrical
statement in which the right hand side of \eqref{eq:11} and indeed
\eqref{eq:12} can be significantly improved -- at a cost!

Kurzweil \cite{MR0073654} showed that, for any non-increasing
function $\psi : \mathbb{N} \rightarrow \mathbb{R}_{>0}$ such that
$\sum \psi(q) = \infty$ and for almost every irrational $\alpha$,
the set of $x$ for which \eqref{eq:1} has infinitely many
solutions is of full Lebesgue measure. This cannot be improved
upon in the sense that there exist irrational $\alpha$ and a
function $\psi$ for which $\sum \psi(q) = \infty$, but the `full
measure' conclusion fails to hold.  Hence, the `almost every'
aspect of Kurzweil's result does not extend to all irrationals
$\alpha$ without modification   -- the divergent sum condition is
not enough.

Over the last few years, there has been much activity in
investigating the shrinking target problem associated with
irrational rotations of the circle. For example, when  $\psi(q) :=
q^{-v}$ ($v > 1$), Bugeaud \cite{MR1972699} and independently
Schmeling $\&$ Trubetskoy \cite{MR1992153} have obtained the
Hausdorff dimension of the set of $x$ for which inequality
\eqref{eq:1} has infinitely many solutions.  Fayad
\cite{MR2268368}, A.-H. Fan $\&$ J.~Wu  \cite{MR2256970}, Kim
\cite{MR2335077} and Tseng \cite{tseng08,tseng_more_remarks} have
built upon the work of Kurzweil in various directions.  In
particular, Kim has proved that for any irrational $\alpha$, the
set of $x$ for which
\begin{equation}
  \label{eq:3}
  \liminf_{q \rightarrow \infty} q \norm{q\alpha-x} = 0
\end{equation}
has full measure. Rather surprisingly, Beresnevich, Bernik, Dodson
$\&$ Velani \cite{beresnevich08:_class_dioph}  have shown that
this result and indeed the  dimension result of Bugeaud and
Schmeling $\&$ Trubetskoy are consequences of the fact that for
any irrational $\alpha$ and any $x$ the inequality \eqref{eq:11}
has infinitely many solutions.

The result of Kim is the underlying motivation for our work. In
this paper we investigate the complementary measure zero set
associated with \eqref{eq:3}; namely
\begin{equation}
  \label{eq:4}
  \Bad_\alpha : = \left\{x \in [0,1) : \exists \text{ } c(x)>0 \text{
      s.t.} \norm{q  \alpha - x} \geq \frac{c(x)}{q} \text{  } \forall
    \text{ } q \in \mathbb{N}\right\}.
\end{equation}

\noindent In fact, we will be concerned with more general actions
than rotations of the circle. Broadly speaking, there are two
natural ways to generalise circle rotations. One option is to
increase the dimension of the torus; i.e. to  consider  the
sequence $\{q \boldsymbol{\alpha}\}$ in $[0,1)^n$ where
$\boldsymbol{\alpha}= (\alpha_1, \dots, \alpha_n)^T \in
\mathbb{R}^n$. The other option is to increase the dimension of
the group acting on the torus; i.e. to consider the sequence
$\{\boldsymbol{\alpha} \cdot \mathbf{q}\}$ where
$\boldsymbol{\alpha}= (\alpha_1, \dots, \alpha_m) \in
\mathbb{R}^m$ and $\mathbf{q} = (q_1, \dots, q_m)^T \in
\mathbb{Z}^m$.

It is possible to consider  both the above mentioned options at
the same time by introducing a $\mathbb{Z}^m$ action on the
$n$-torus by $n \times m$ matrices. Indeed, we may consider the
points $\{A\mathbf{q}\} \in [0,1)^n$ where $A \in \mat_{n \times
  m}(\mathbb{R})$ is fixed and $\mathbf{q}$ runs over $\mathbb{Z}^m$.
In this case, the natural analogue of  $\Bad_\alpha$  is the set
\begin{equation*}
  \label{eq:5}
  \Bad_A : = \left\{\mathbf{x} \in [0,1)^n : \exists \text{ }
    c(\mathbf{x})>0 \text{ s.t.} \norm{A\mathbf{q} -
      \mathbf{x}} \geq \frac{c(\mathbf{x})}{\infabs{\mathbf{q}}^{m/n}}
    \text{  } \forall \text{ }\mathbf{q} \in \mathbb{Z}^m \setminus \left\{ \bf{0} \right\} \, \right\}.
\end{equation*}

\noindent Here and throughout, for a vector $\mathbf{x}$ in
$\mathbb{R}^n$ we will denote by $\infabs{\mathbf{x}}$ the maximum
of the absolute values of the coordinates of $\mathbf{x}$; i.e.
the infinity norm of $\mathbf{x}$. Also, $\norm{\mathbf{x}} :=
\min_{\mathbf{y} \in \mathbb{Z}^n}\infabs{\mathbf{x}-\mathbf{y}}$.

\vspace{1ex}

The underlying goal of this paper is to show that no matter which
of the $\mathbb{Z}^m$ actions defined above we choose, the set
$\Bad_A$ is of maximal Hausdorff dimension.

\begin{thm} \label{startthm}
For any  $A \in \mat_{n \times m}(\mathbb{R})$,  $$ \dim \Bad_A \,
= \, n \, . $$
\end{thm}

\noindent In terms of the more familiar setting of irrational
rotations of the circle, the theorem  reads as follows.

\begin{cor}
  \label{cor:1Dbad}
  For any $\alpha \in \mathbb{R}$,
  \begin{equation*}
    \hdim \, \Bad_\alpha = 1 \, .
  \end{equation*}
\end{cor}


\noindent Note that if  $\alpha$ is rational, the set
$\Bad_\alpha$ is easily seen to contain all points in the unit
interval bounded away from a finite set of points. Thus, for
rational $\alpha$ not only is $\Bad_\alpha$ of full dimension but
it is of full Lebesgue measure. In higher dimensions, similar
phenomena occur in which the finite set of points is replaced by a
finite set of affine subspaces. The reader is referred to
\cite{MR2266457} and \S\ref{sec:construction-of-sequences} below
for  further details.



\vspace{1ex}

 Inspired by the works of  Kleinbock $\&$ Weiss
\cite{MR2191212} and Kristensen, Thorn $\&$ Velani
\cite{MR2231044},  we shall deduce  Theorem \ref{startthm} as a
simple consequence of a general statement concerning the
intersection of $\Bad_A$ with compact subsets of $\mathbb{R}^n$.
The latter includes exotic fractal sets such as the Sierpinski
gasket and the van Koch curve.

\vspace{1ex}

\section{The setup and main result}
\label{sec:statement-results}


  Let $(X, d)$ be
a metric space and $(\Omega, d)$ be a compact subspace of $X$
which supports a non-atomic finite measure $\mu$. Throughout,
$B(c, r)$ will denote a closed ball in $X$ with center $c$ and
radius $r$. The measure $\mu$ is said to \emph{$\delta$-Ahlfors
regular} if there exist strictly positive constants $\delta$ and
$r_0$ such that for $c \in \Omega$ and $r< r_0$
\begin{equation*}
 a r^{\delta} \leq \mu(B(c, r)) \leq b r^{\delta}   \ ,
\end{equation*}
where $0 < a \leq 1 \leq b$ are constants independent of the ball.
It is easily verified that if $\mu$ is $\delta$-Ahlfors regular
then the Hausdorff dimension of $\Omega$ is $\delta$; i.e.
\begin{equation}
\label{regcon} \dim \Omega = \delta  \ .
\end{equation}
For  further details including  the definition of Hausdorff
dimension the reader is referred to  \cite{Mat}.

In the above, take $X = \mathbb{R}^n$ and let $\mathcal{L}$ denote
a generic $(n-1)$-dimensional hyperplane. For $\epsilon > 0$, let
$\mathcal{L}^{(\epsilon )}$ denote the $\epsilon$-neighbourhood of
$\mathcal{L}$.  The measure $\mu$ is said to be  \emph{absolutely
$\alpha$-decaying} if there exist
   strictly  positive constants $C, \alpha$ and $r_0$ such that for any
    hyperplane $\mathcal{L}$, any $\epsilon >0$, any $x \in \Omega$
    and any $r< r_0$,
    \begin{equation*}
      \mu(B(x, r) \cap \mathcal{L}^{(\epsilon )}) \leq C \left(
        \frac{\epsilon}{r} \right)^{\alpha}
      \mu(B(x, r)).
    \end{equation*}

\noindent 
It is worth mentioning that if $\mu $ is $\delta$-Ahlfors regular
and absolutely $\alpha$-decaying, then $\mu$ is an absolutely
friendly measure as defined in \cite{PV}.

Armed with the notions of Ahlfors regular and absolutely decaying,
we are in the position to state our main result.

\vspace*{1ex}

\begin{thm}
  \label{thm:nDbad}
  Let $K \subseteq [0,1]^n$ be a compact set which supports an absolutely
  $\alpha$-decaying, $\delta$-Ahlfors regular measure $\mu$  such that
  $\delta > n-1$.  Then, for any $A \in \mat_{n \times
    m}(\mathbb{R})$,
  \begin{equation*}
    \hdim (\Bad_A \cap K) =  \delta.
  \end{equation*}
\end{thm}

\vspace*{1ex}

\noindent In view of (\ref{regcon}), the theorem  can be
interpreted as stating that within $K$ the set $\Bad_A$  is of
maximal dimension.
 It is easily seen that  Theorem \ref{startthm} is a consequence of Theorem \ref{thm:nDbad} -- simply take
$K= [0,1]^n$ and $\mu$ to be $n$-dimensional Lebesgue measure.
Trivially,  $n$-dimensional Lebesgue measure is $n$-Ahlfors
regular and absolutely $1$-decaying. More exotically, the natural
measures associated with self-similar sets in $ \mathbb{R}^n$
satisfying the open set condition are absolutely $\alpha$-decaying
and $\delta$-Ahlfors regular -- see  \cite{MR2191212,PV}. Thus,
Theorem \ref{thm:nDbad} is applicable to these sets which in
general are of fractal nature.

\vspace*{2ex}

Although Theorem \ref{thm:nDbad} constitutes our main result, we
state an `auxiliary' result in this section for the simple fact
that it is new and of independent interest. In short, it strengthens and
generalises a  theorem  of Pollington \cite{MR540398} and de
Mathan \cite{MR510723, MR612195} that  answers a question of
Erd\H{o}s. 
A sequence $ \textstyle { \left\{ \mathbf{y}_i \right\} := \left\{
\mathbf{y}_i := (y_{1,i}, \dots, y_{n,i})^T \in \mathbb{Z}^n
\setminus \left\{\bf{0}\right\} \ \right\} }  $
is said to be \emph{lacunary} if there exits a constant $
\lambda>1 $ such that
$$
\infabs{\mathbf{y}_{i+1}}  \; \geq \; \lambda \,
\infabs{\mathbf{y}_i}
 \qquad \forall \ i \in \mathbb{N} \ .
$$
%
  Given a sequence  $\left\{ \mathbf{y}_i \right\}$ in $\mathbb{Z}^n
$, let
\begin{equation*}
    \Bad_{\left\{ \mathbf{y}_i \right\} } :=\left\{ \mathbf{x} \in [0,1]^n :
      \exists \text{ } c(\mathbf{x})>0 \, \text{ s.t.}  \, \norm{\mathbf{y}_i
        \cdot \mathbf{x}} \geq c(\mathbf{x})  \text{ } \forall \text{
      } i \in \mathbb{N}\right\}  \ .
  \end{equation*}

\vspace*{1ex}

\begin{thm}
  \label{lem:pollington}
  Let $\left\{ \mathbf{y}_i \right\} $ be a lacunary sequence  in
  $\mathbb{Z}^n$.
  Furthermore, let $K \subseteq [0,1]^n$ be a
  compact set which supports an absolutely $\alpha$-decaying,
  $\delta$-Ahlfors regular measure $\mu$ such that $\delta > n-1$.
  Then
  \begin{equation*}
    \hdim (\Bad_{\left\{\mathbf{y}_i\right\}} \cap K )  = \delta.
  \end{equation*}
\end{thm}

\vspace{1ex}

\noindent On setting $n=1$, $K = [0,1]$ and  $\mu$ to be
one-dimensional Lebesgue measure, Theorem \ref{lem:pollington}
corresponds to the theorem of Pollington and de Mathan referred to
above.

\section{Preliminaries for Theorem \ref{lem:pollington} }
\label{sec:framework}

The proof of Theorem \ref{lem:pollington} makes use of the general
framework developed in \cite{MR2231044} for establishing dimension
statements for a large class of badly approximable sets. In this
section we provide a simplification of the framework that is
geared towards the particular application we have in mind. In
turn, this will avoid excessive referencing to the conditions
imposed in \cite{MR2231044} and thereby improve the clarity of our
exposition.

As in \S\ref{sec:statement-results}, let $(X, d)$ be a metric
space and $(\Omega, d)$ be a compact subspace of $X$ which
supports a non-atomic finite measure $\mu$. Let $\mathcal{R} :=
\left\{ R_{\alpha} \in X : \alpha \in J \right\} $ be a family of
subsets $R_{\alpha}$ of $X$ indexed by an infinite countable set
$J$.  The sets $R_{\alpha}$ will be referred to as the
\emph{resonant sets}. Next, let $\beta :J\rightarrow
\mathbb{R}_{>0}:\alpha \mapsto \beta_{\alpha}$ be a positive
function on $J$ such that the number of $\alpha \in J$ with $
\beta_{\alpha} $ bounded above is finite. Thus, $ \beta_{\alpha} $
tends to infinity as $\alpha$ runs through $J$. We are now in the
position to define the badly approximable set
\begin{equation*}
  \Bad(\mathcal{R}, \beta) := \left\{x \in \Omega : \exists \text{
    } c(x)>0 \text{ s.t. } d(x, R_{\alpha} ) \geq
    \frac{c(x)}{\beta_{\alpha}} \text {  } \forall \text{ } \alpha \in
    J \right\}  \ ,
\end{equation*}
where   $d(x, R_{\alpha} ):= \inf_{a \in R_{\alpha}}d(x, a)$.
Loosely speaking, $\Bad(\mathcal{R}, \beta)$ consists of points in
$\Omega $ that `stay clear' of the family $ \mathcal{R} $ of
resonant sets by a factor governed by $\beta$.

The goal is to determine conditions under which $\dim
\Bad(\mathcal{R}, \beta)  =  \dim \Omega $; that is to say that
the set of badly approximable points in $\Omega$ is of  maximal
dimension.  With this in mind, we begin with some useful notation.
For any fixed integer $k>1$ and any integer  $n \geq 1$, let $B_n
:= \left\{ x \in \Omega : d(c, x)
  \leq 1/k^n \right\}$ denote a generic closed ball in $\Omega$ of radius $1/k^n$
with centre $c$ in $\Omega$.  For any $\theta \in
\mathbb{R}_{>0}$, let $\theta B_n := \left\{ x \in \Omega : d(c,
x)
  \leq \theta/k^n \right\}$ denote the ball $B_n$ scaled by $\theta$.
Finally,  let $J(n) := \left\{ \alpha \in J : k^{n-1} \leq
\beta_{\alpha} < k^n \right\}$. The following statement is a
simple consequence of  combining Theorem 1 and Lemma 7 of
\cite{MR2231044} and realises the above mentioned goal.

\begin{thmktv}
  \label{thm:1fromthat}
  Let $(X, d)$ be
a metric space and $(\Omega, d)$ be a compact subspace of $X$
which supports of a $\delta$-Ahlfors
  regular measure $\mu$. Let $k$ be sufficiently large. Then for any
  $\theta \in \mathbb{R}_{>0}$, any $ n \geq 1 $ and any ball $B_n $
  there exists a collection $\mathcal{C}(\theta B_n)$ of
  disjoint balls $2 \theta B_{n+1}$ contained within $\theta B_n$ such
  that $ \# \mathcal{C}(\theta B_n) \geq \kappa_1 \,  k^{\delta} $  .
 In addition, suppose for  some  $\theta \in \mathbb{R}_{>0}$ we also have that
  \begin{equation}
    \label{eq:cond2}
    \# \left\{ 2 \theta B_{n+1} \subset \mathcal{C}(\theta B_n):
      \min_{\alpha \in J(n+1)} d(c, R_{\alpha}) \leq 2
        \theta k^{-(n+1) } \right\} \leq \kappa_2 k^{\delta}  \, ,
  \end{equation}
  where $0< \kappa_2 < \kappa_1$  are absolutely constants independent of $k$ and $n$.  Furthermore, suppose
  \begin{equation}
    \label{eq:cond3}
    \hdim \, \left( \cup_{\alpha \in J} R_{\alpha} \right) < \delta  \, .
  \end{equation}
  Then
  \begin{equation*}
    \hdim \, \Bad(\mathcal{R}, \beta) = \delta  \, .
  \end{equation*}
\end{thmktv}

\vspace*{2ex}

\noindent Note that the theorem together with  (\ref{regcon})
implies that  $\dim \Bad(\mathcal{R}, \beta)  =  \dim \Omega $.

\vspace*{2ex}

\section{Proof of Theorem \ref{lem:pollington}}
\label{sec:proof-of-lemma-pollington}

  We are given a lacunary sequence $\left\{ \mathbf{y}_i \right\}$. For
 each index $i \in \mathbb{N} $ and any integer $p$,  consider the
   hyperplane $\mathcal{L}_{p, i}:= \{\mathbf{x}  \in \mathbb{R}^n : \mathbf{y}_i \cdot \mathbf{x} = p\}$ .
  It is easily verified that
  $\Bad_{\left\{\mathbf{y}_i\right\}} \cap K $ is equivalent to the set of
  $\mathbf{x}$ in  $K$ for which there exists a constant
  $c(\mathbf{x})>0$ such that $\mathbf{x}$ avoids the
  $c(\mathbf{x})/\infeuc{\mathbf{y}_i}$--neighbourhood of
  $\mathcal{L}_{p, i}$ for every choice of $i$ and $p $; that is
  \begin{equation*}
    \Bad_{\left\{\mathbf{y}_i\right\}} \cap K  = \left\{ \mathbf{x} \in K :
      \exists  \text{ } c(\mathbf{x})>0 \, \text{ s.t.} \,
      \min_{\mathbf{y}\in \mathcal{L}_{p, i}}
      \infeuc{\mathbf{x}-\mathbf{y}} \geq
      \frac{c(\mathbf{x})}{\infeuc{\mathbf{y}_i}}
      \text{ } \forall \text{ }(p, i) \in
      \mathbb{Z}\times \mathbb{N} \right\}.
  \end{equation*}

\noindent Here $\infeuc{\, . \, }$  is the standard Euclidean norm
in $\mathbb{R}^n$. With reference to \S\ref{sec:framework},  set
\begin{eqnarray*}
& &X:=\mathbb{R}^n \, ,\quad \Omega :=K \, ,\quad  d:= \infeuc{\,
. \, } \ , \quad J:=\left\{(p, i) \in
      \mathbb{Z}\times \mathbb{N} \right\} \, ,
\\ 
& & \alpha:=(p, i) \in J \, , \quad
         R_{\alpha}:= \mathcal{L}_{p, i} \quad
    \rm{and} \quad    \beta_{\alpha}:=\infeuc{\mathbf{y}_i}  \ .
\end{eqnarray*}

\noindent  It follows that
  $$\Bad(\mathcal{R}, \beta) = \Bad_{\left\{\mathbf{y}_i\right\}} \cap K  \ . $$

\noindent The upshot of this is that the proof of Theorem
\ref{lem:pollington} is reduced to  showing that the conditions of
Theorem KTV are satisfied.

For $k > 1$ and $m \geq 1$, let  $B_m $ be a generic closed ball
of radius $k^{-m}$ and centre in $K$.  For $k $ sufficiently large
and any $\theta \in \mathbb{R}_{>0}$, Theorem KTV guarantees the
existence of  a collection $\mathcal{C}(\theta B_m)$ of disjoint
balls $2 \theta B_{m+1}$ contained within $\theta B_m$ such that
\begin{equation*}
    \label{eq:cond1}
    \# \mathcal{C}(\theta B_m) \geq \kappa_1 \, k^{\delta} \,  .
\end{equation*}
The positive constant $\kappa_1$ is independent of $k$ and $n$. We
now endeavor to show that the additional condition
\eqref{eq:cond2} on the collection $\mathcal{C}(\theta B_m)$ is
satisfied.  To this end,  set $\theta: = (2 k )^{-1}$ and  proceed
as follows.
  Fix $m \geq 1$ and assume that there exists
  an index $i$ such that
  \begin{equation}
  \label{help}
  k^m \leq \infeuc{\mathbf{y}_i} < k^{m+1}  \, .
\end{equation}
If this was not the case,  the left hand side of \eqref{eq:cond2}
is zero and the additional condition  is  trivially  satisfied.
Associated with the index $i$ is the family of hyperplanes
$\left\{ \mathcal{L}_{p, i} : p \in   \mathbb{Z} \right\} $. The
distance between any two such hyperplanes  is at least
$\infeuc{\mathbf{y}_i}^{-1} > k^{-(m+1)}$.  The diameter of the
ball $\theta B_m$ is $ k^{-(m+1)}$. Thus, for any element of the
sequence  $\left\{ \mathbf{y}_i \right\}$ satisfying (\ref{help})
there is at most one hyperplane passing through  $\theta B_m$.
Assume, the hyperplane $\mathcal{L}_{p, i}$ passes through
$\theta B_m$ and consider the counting function
  \begin{equation*}
    \omega (m, p, i):=\# \left\{2\theta B^{m+1} \subset
    \mathcal{C}(\theta B_m) : 2\theta B_{m+1} \cap \mathcal{L}_{p, i}
    \neq \emptyset\right\}   \, .
  \end{equation*}
The balls  $2 \theta B_{m+1}$ are disjoint and each is of diameter
$4\theta k^{-(m+1)}$.  Thus,  on setting $\epsilon: =8\theta
k^{-(m+1)}$  it follows that
 \begin{eqnarray*}
    \omega (m, p, i) & \leq &  \# \left\{2\theta B_{m+1} \subset
      \mathcal{C}(\theta B_m) : 2\theta B_{m+1} \subset
      \mathcal{L}^{(\epsilon )}_{p, i} \right\}  \\
      & \leq &
  \frac{\mu(\theta B_m \cap \mathcal{L}^{(\epsilon )}_{p,
        i})}{\mu(2\theta B_{m+1})}    \ .
 \end{eqnarray*}
On making use of  the fact that $\mu$ is absolutely
$\alpha$-decaying and $\delta$-Ahlfors regular, it is readily
verified that
\begin{equation*}
    \omega (m, p, i) \,  \leq  \, \kappa  \, k^{\delta - \alpha}   \  .
\end{equation*}
The absolute constant  $\kappa$ is  dependent only on  $\alpha$
and $\delta$.  Next, let  $\upsilon(m,\left\{ \mathbf{y}_i
\right\})   $ denote the number of elements of the sequence
$\left\{ \mathbf{y}_i \right\}$  satisfying (\ref{help}).
Since  $\left\{ \mathbf{y}_i \right\}$ is lacunary, we find that
for $k$ sufficiently large
$$
\upsilon(m,\left\{ \mathbf{y}_i \right\})  \leq  1+ \log (\sqrt{n}
\, k)/\log \lambda  \, < \, \frac{\kappa_1}{2 \, \kappa}
k^{\alpha}  \, . $$ Here, $\lambda > 1 $ is the lacuarity constant
and we have used the fact that $  \infabs{\mathbf{y}} \leq
\infeuc{\mathbf{y}}  \leq \sqrt{n} \, \infabs{\mathbf{y}}$ for $
\mathbf{y} \in \mathbb{Z}^n $.  On combining the above upper bound
estimates,  we have that
\begin{eqnarray*}
  \text{l.h.s. of \eqref{eq:cond2}}  & < &   \upsilon(m,\left\{ \mathbf{y}_i \right\})    \ \times \ \omega (m, p, i)  \\
      & \leq & \frac{\kappa_1}{2 \kappa} k^{\alpha}
   \ \times \ \kappa k^{\delta - \alpha}  \ = \   \textstyle{\frac12} \,  \kappa_1  k^{\delta} \ .
 \end{eqnarray*}
 Thus, with $\theta: = (2 k )^{-1}$  the collection $\mathcal{C}(\theta B_m)$ satisfies  \eqref{eq:cond2}. Finally, note that  the  collection $\left\{  \mathcal{L}_{p, i} : (p, i) \in \mathbb{Z} \times \mathbb{N} \right\}$ of hyperplanes (resonant sets)  is countable  and so
  \begin{equation*}
    \hdim \left( \cup \,  \mathcal{L}_{p, i} \right) \, = \, n-1    \ .
  \end{equation*}
We  are  given that $\delta > n-1$ and so  \eqref{eq:cond3} is
trivially satisfied.   Thus, the conditions of Theorem KTV are
satisfied and Theorem  \ref{lem:pollington} follows.

\section{Preliminaries for Theorem \ref{thm:nDbad}}
\label{sec:construction-of-sequences}

The proof of Theorem \ref{thm:nDbad} makes use of the existence of
`special' sequences which for the most part are constructed in
 \cite{MR2266457}. Throughout, $\mat_{n \times
m}^*(\mathbb{R})$  will denote the collection of matrices  $A \in
\mat_{n \times m}(\mathbb{R})$ for which
the associated group $G := A^T \mathbb{Z}^n + \mathbb{Z}^m$ has
rank $n+m$. In Section 3 of \cite{MR2266457}, it is shown that
associated with each matrix $ A \in \mat_{n \times
m}^*(\mathbb{R})$ there exists a sequence $\left\{ \mathbf{y}_i
\right\}$ of integer vectors $\mathbf{y}_i = (y_{1,i}, \dots,
y_{n,i})^T  \in \mathbb{Z}^n$ satisfying the following properties:

\begin{itemize}
\item[(i)] $1 = \infabs{\mathbf{y}_1} < \infabs{\mathbf{y}_2} <
  \infabs{\mathbf{y}_3} < \dots$ \, ,

\vspace{2ex}

  \item[(ii)] $\norm{A^T \mathbf{y}_1} > \norm{A^T \mathbf{y}_2} >
    \norm{A^T \mathbf{y}_3} > \dots$ \, ,

\vspace{2ex}

  \item[(iii)] For all non-zero $\mathbf{y} \in \mathbb{Z}^n$ with
    $\infabs{\mathbf{y}} < \infabs{\mathbf{y}_{i+1}}$ we have that
    $\norm{A^T \mathbf{y}} \geq \norm{A^T \mathbf{y}_i}$  .
\end{itemize}

\noindent Such a sequence $\left\{ \mathbf{y}_i \right\}$  is
referred to as a \emph{sequence of best approximations} to $A$. In
the one-dimensional case ($n=m=1$), when $A $ is an irrational
number $\alpha$,  the sequence of best approximations is precisely
the sequence of denominators associated with the convergents of
the continued fraction representing $\alpha$.

\noindent Let $\left\{ \mathbf{y}_{i} \right\}$ be a sequence of
best approximations to a matrix $A \in \mat_{n \times m}^*
(\mathbb{R})$. A further property enjoyed by $\left\{ \mathbf{y}_i
\right\}$, is that
\begin{equation}
  \label{eq:mink}
  \norm{A^T \mathbf{y}_i} \leq \infabs{\mathbf{y}_{i+1}}^{-m/n} \;
  \text{  }  \forall \text{ } i \in \mathbb{N}.
\end{equation}
\noindent This property is easily deduced via  Dirichlet's box
principle  -- see Section 3 of \cite{MR2266457} for the details.

\vspace*{1ex}

The following result, which is taken from Section 5 of
\cite{MR2266457},  enables us to extract a lacunary subsequence
from a given sequence of best approximations. This  will allow us
to utilise  Theorem \ref{lem:pollington} in the course of
establishing  Theorem \ref{thm:nDbad}.

\begin{lembl}
\label{lem:lacsubseq} Let $A \in \mat_{n \times m}^*(\mathbb{R})$
and  let $\left\{ \mathbf{y}_{i} \right\} $ be a sequence of best
approximations to  $A$.  Then, there exists an increasing function
$\phi : \mathbb{N}\rightarrow \mathbb{N}$ such that  $\phi (1) =
1$ and for  $i \geq 2$
\begin{equation}
\label{eq:8}
    \infabs{\mathbf{y}_{\phi(i)}} \geq \sqrt{9n} \,
    \infabs{\mathbf{y}_{\phi(i-1)}} \quad
    \text{and} \quad   \infabs{\mathbf{y}_{\phi(i-1)+1}} \geq
    \frac{ \infabs{\mathbf{y}_{\phi(i)}}}{9n }  \ .
\end{equation}
\end{lembl}

\vspace{1ex} \noindent It is clear that the sequence $\left\{
\mathbf{y}_{\phi(i)} \right\}$ is lacunary and that it also
satisfies \eqref{eq:mink}; i.e.
\begin{equation}
\label{eq:7}
  \norm{A^T \mathbf{y}_{\phi(i)}} \leq \infabs{\mathbf{y}_{\phi(i)+1}}^{-m/n} \;
 \text{  }  \forall \text{ } i \in \mathbb{N}.
\end{equation}

The next inequality   follows  directly from the definition of the
norms involved. For any $\mathbf{x}$ and $\mathbf{y}$ in
$\mathbb{R}^k$, we have that
\begin{equation}
  \label{eq:norm}
  \norm{\mathbf{x} \cdot \mathbf{y}}  \,  < \,  k  \,
  \infabs{\mathbf{x}}\norm{\mathbf{y}}  \ .
\end{equation}

We end this section  with a short discussion that allows us to
assume that $A \in \mat_{n \times m}^*(\mathbb{R})$ when proving
Theorem \ref{thm:nDbad}. With this in mind, suppose $A \in \mat_{n
\times m}(\mathbb{R})$ and that the rank of the associated group
$G := A^T \mathbb{Z}^n + \mathbb{Z}^m$ is strictly less than
$n+m$.  Then, it is easily verified that $\left\{ A \mathbf{q}:
\mathbf{q} \in \mathbb{Z}^m\right\}$ is restricted to at most a
countable  family of positively separated,  parallel hyperplanes
in $\mathbb{R}^n$. Let $ F$ denote the set of these hyperplanes.
Then,
$$
K \setminus F = \Bad_A \cap K   \ .
$$
We are given that $\delta
> n-1$ which together with (\ref{regcon}) implies that $\dim K $ is strictly
greater than $\dim F$.  Thus,  $\dim (K \setminus F) = \dim K $
and the statement of Theorem \ref{thm:nDbad} follows for any $A
\notin \mat_{n \times m}^*(\mathbb{R})$.

\section{Proof of Theorem \ref{thm:nDbad}}
\label{sec:proof-of-theorem-1}

Without loss of generality, assume that $A \in \mat_{n \times
m}^*(\mathbb{R})$ and let $\left\{ \mathbf{y}_i \right\}$ be a
sequence of best approximations to $A$. In view of Lemma BL, there
exists a lacunary subsequence $\left\{ \mathbf{y}_{\phi(i)}
\right\}$ of the sequence of best approximations. For any $c > 0$,
let
\begin{equation*}
    {\bf B}_{\left\{ \mathbf{y}_{\phi(i)}
\right\}}(c): = \left\{\mathbf{x} \in K \colon
    \norm{\mathbf{y}_{\phi(i)} \cdot \mathbf{x}} \geq c  \text{ } \forall \text{ } i \in \mathbb{N}\right\}  \, .
\end{equation*}
It is readily verified that $ \Bad_{\left\{ \mathbf{y}_{\phi(i)}
\right\}} \cap K  \, =   \, \bigcup_{c > 0}   {\bf B}_{\left\{
\mathbf{y}_{\phi(i)} \right\}}(c) $
and that
\begin{equation*}
\dim   {\bf B}_{\left\{ \mathbf{y}_{\phi(i)} \right\}}(c) \ \to  \
\dim \, ( \Bad_{\left\{ \mathbf{y}_{\phi(i)} \right\}} \cap K )
\quad \text{ as  }  \  c \to 0  \, .
\end{equation*}

For $c$ sufficiently small, suppose for the moment that
\begin{equation}
\label{assup}
  {\bf B}_{\left\{ \mathbf{y}_{\phi(i)}
\right\}}(c) \ \subseteq \  \Bad_A \cap K  \ .
\end{equation}
On utilising Theorem \ref{lem:pollington}, it follows that
\begin{equation*}
\dim \, (\Bad_{\left\{ \mathbf{y}_{\phi(i)} \right\}}  \cap K ) \,
\geq  \, \dim {\bf B}_{\left\{ \mathbf{y}_{\phi(i)} \right\}}(c)
\, \to \, \delta  \quad \text{ as  }  \  c \to 0  \, .
\end{equation*}
The upshot of this is that $ \dim \, (\Bad_{\left\{
\mathbf{y}_{\phi(i)} \right\}}  \cap K ) \, \geq \, \delta $. For
the complementary upper bound statement,  trivially
$$
\dim \, (\Bad_{\left\{ \mathbf{y}_{\phi(i)} \right\}}  \cap K ) \,
\leq \, \dim \, K \, \stackrel{(\ref{regcon})}{=} \, \delta \, .
$$
This completes the proof of Theorem \ref{thm:nDbad} modulo the
inclusion (\ref{assup}).

To establish (\ref{assup}), fix a point $\mathbf{x}$ in ${\bf
B}_{\left\{ \mathbf{y}_{\phi(i)} \right\}}(c)$ and let $\mathbf{q}
$ be any non-zero integer vector. For $c$ sufficiently small,
there exists an index $i \in \mathbb{N}$ such that
\begin{equation}
    \label{eq:9}
    \infabs{\mathbf{y}_{\phi(i)}} \, \leq \,  9 n
    \left(\frac{2m}{c}\right)^{m/n} \infabs{\mathbf{q}}^{m/n}   \, < \,
    \infabs{\mathbf{y}_{\phi(i + 1)}}.
\end{equation}
The existence of such an index is guaranteed by the first of the
inequalities in \eqref{eq:8} as long as  $c$ is sufficiently
small. By the  definition of ${\bf B}_{\left\{\mathbf{y}_{\phi(i)}
\right\}}(c)$ and the trivial equality
\begin{equation*}
    \mathbf{y}_{\phi(i)} \cdot \mathbf{x} \, = \, \mathbf{q} \cdot A^T
    \mathbf{y}_{\phi(i)} - \mathbf{y}_{\phi(i)}\cdot (A \mathbf{q} -
    \mathbf{x}),
\end{equation*}
we immediately have that
\begin{equation}
  \label{eq:triv}
  0 < c \leq \norm{\mathbf{y}_{\phi(i)} \cdot \mathbf{x}} = \norm{\mathbf{q} \cdot A^T
    \mathbf{y}_{\phi(i)} - \mathbf{y}_{\phi(i)}\cdot (A \mathbf{q} -
    \mathbf{x})} \, .
\end{equation}
On applying the triangle inequality and making use of
\eqref{eq:norm}, it follows that
  \begin{equation}
    \label{eq:2}
    c \leq m \infabs{\mathbf{q}} \norm{A^T \mathbf{y}_{\phi(i)}} + n
    \infabs{\mathbf{y}_{\phi(i)}} \norm{A \mathbf{q} - \mathbf{x}}.
  \end{equation}
However,
\begin{equation*}
    \label{eq:10}
    m \infabs{\mathbf{q}} \norm{A^T \mathbf{y}_{\phi(i)}} \, \stackrel{(\ref{eq:7})}{\leq} \,   m
    \infabs{\mathbf{q}} \infabs{\mathbf{y}_{\phi(i)+1}}^{-n/m} \, \stackrel{(\ref{eq:9})}{\leq} \,
    \frac{m}{(9 n)^{n/m} \frac{2m}{c}}
    \left(\frac{ \infabs{\mathbf{y}_{\phi(i+1)}} }{
    \infabs{\mathbf{y}_{\phi(i)+1}}}\right)^{n/m}   \, \stackrel{(\ref{eq:8})}{\leq} \,  \frac{c}{2}
\end{equation*}
and
  \begin{equation*}
    \label{eq:6}
    n \infabs{\mathbf{y}_{\phi(i)}} \norm{A \mathbf{q} -\mathbf{x}}
   \,  \stackrel{(\ref{eq:9})}{\leq}  \,  9 n^2 \left(\frac{2m}{c}\right)^{m/n}
    \infabs{\mathbf{q}}^{m/n} \norm{A \mathbf{q}-\mathbf{x}}   \ ,
  \end{equation*}
which together with \eqref{eq:2} yields  that
  \begin{equation*}
    \norm{A \mathbf{q} -\mathbf{x}} >
    \frac{c^{m/n+1}}{9n^2(2m)^{m/n}} \infabs{\mathbf{q}}^{-m/n}.
  \end{equation*}
In other words, for any $c$ sufficiently small
\begin{equation*}
    {\bf B}_{\left\{ \mathbf{y}_{\phi(i)}
\right\}}(c) \subseteq \left\{\mathbf{x} \in K : \exists  \text{ }
      c(\mathbf{x})>0 \text{ s.t.} \norm{A\mathbf{q} -
        \mathbf{x}} \geq
      \frac{c(\mathbf{x})}{\infabs{\mathbf{q}}^{m/n}} \text{ } \forall \text{ }
      \mathbf{q} \in \mathbb{Z}^m \setminus \left\{ \bf{0} \right\}    \right\}.
\end{equation*}
The right hand side is $\Bad_A \cap K$ and this establishes
(\ref{assup}) which in turn completes the proof of Theorem
\ref{thm:nDbad}.

\providecommand{\bysame}{\leavevmode\hbox
to3em{\hrulefill}\thinspace}
\providecommand{\MR}{\relax\ifhmode\unskip\space\fi MR }
\providecommand{\MRhref}[2]{%
  \href{http://www.ams.org/mathscinet-getitem?mr=#1}{#2}
} \providecommand{\href}[2]{#2}

\end{document}